\documentclass[12pt,a4paper]{amsart}
\usepackage[all]{xy}
\usepackage{amsthm}
\newtheorem{THM}{Theorem}

\newtheorem{thm}{Theorem}
\newtheorem{lem}[thm]{Lemma}
\newtheorem{prop}[thm]{Proposition}
\newtheorem{cor}[thm]{Corollary}
\theoremstyle{definition}

\theoremstyle{remark}

\newtheorem{ex}[thm]{Example}
\DeclareMathAlphabet{\mathbfit}{OT1}{cmr}{bx}{it}
\newcommand{\iso}{\cong}
\newcommand{\Z}{\mathbf{Z}}
\newcommand{\Q}{\mathbf{Q}}

\newcommand{\F}[1]{\mathbf{F}_{\!{#1}}}
\newcommand{\Zloc}[1]{\Z_{(#1)}}

\newcommand{\Lie}[1]{L_{#1}}
\newcommand{\qi}{\xrightarrow{\simeq}}
\newcommand{\dual}{\sharp}
\newcommand{\torsion}{\operatorname{Torsion}}
\newcommand{\free}[1]{F{#1}}
\newcommand{\image}{\operatorname{im}}
\newcommand{\Hom}{\operatorname{Hom}}
\newcommand{\Der}{\operatorname{Der}}
\newcommand{\rank}{\operatorname{rank}}
\newcommand{\cw}[2]{\mathbf{CW}^{#1}_{\!{#2}}}

\begin{document}
\title{A Torsion-Free Milnor-Moore Theorem}
\author{Jonathan A. Scott}
\address{Department of Mathematical Sciences \\
	University of Aberdeen \\
	Aberdeen AB24 3UE \\
	United Kingdom}
\email{j.scott@maths.abdn.ac.uk}
\thanks{This article was written while the author was an NSERC 
	Postdoctoral Fellow.}
\date{\today}

\begin{abstract}
Let $\Omega X$ be the space of Moore loops on a finite, $q$-connected, 
$n$-dimensional CW complex $X$, and let $R\subset\Q$ be a subring containing
$1/2$.
Let $\rho(R)$ be the least non-invertible prime in $R$.
For a graded $R$-module $M$ of finite type, let 
$\free{M} = M / \torsion{M}$.
We show that the inclusion 
$P \subset \free{H_{*}(\Omega X;R)}$ of the sub-Lie algebra of
primitive elements induces
an isomorphism of Hopf algebras
$UP \xrightarrow{\iso} \free{H_{*}(\Omega X;R)}$,
provided $\rho(R) \geq n/q$.
Furthermore, the Hurewicz homomorphism induces an
embedding of $F(\pi_{*}(\Omega X)\otimes R)$ in $P$,
with $P / F(\pi_{*}(\Omega X)\otimes R)$ torsion.
As a corollary, if $X$ is elliptic, then $\free{H_{*}(\Omega X;R)}$
is a finitely generated $R$-algebra.
\end{abstract}

\maketitle


\section{Introduction}\label{sec:introduction}

Let $\Omega X$ be the Moore loop space on a pointed, simply-connected
topological space $X$.
A theorem of Milnor and Moore~\cite{milnor-moore:65}
states that $H_{*}(\Omega X ; \Q)$ is the universal enveloping algebra
of its sub-Lie algebra of primitive
elements.
Halperin~\cite{halperin:92} established the same result working over
any subring $R \subset \Q$ containing $1/2$, provided $X$ is a finite
CW complex, $H_{*}(\Omega X; R)$ is torsion-free, and the least
non-invertible prime in $R$ is sufficiently large.
For a graded $R$-module $M$ of finite type, set 
$\free{M} = M / \torsion{M}$.
Let $\cw{n}{q}$ be the full subcategory of pointed topological
spaces consisting of all finite, $q$-connected ($q \geq 1$),
$n$-dimensional CW complexes.
Suppose $X \in \cw{n}{q}$, and
let $P$ be the sub-Lie algebra of primitive elements of the Hopf algebra
$\free{H_{*}(\Omega X;R)}$.
Denote by $\rho(R)$ the least non-invertible prime in $R$.
The aim of this paper is to prove the following generalisation of Halperin's
result.
 
\begin{THM}\label{thm:bigyin}
With the hypotheses and notation above,
if $\rho(R) \geq n/q$, then the inclusion 
$P \subset \free{H_{*}(\Omega X;R)}$
extends to an isomorphism of Hopf algebras
$UP \xrightarrow{\iso} \free{H_{*}(\Omega X;R)}$.
\end{THM}

Theorem~\ref{thm:bigyin} should not be surprising.
Popescu~\cite{popescu:00} proved it for two-cones
$X$ with $\dim{X}<\rho(R)$.  
The present author showed that if $X \in \cw{n}{q}$ and $p \geq n/q$,
then $\free{H_{*}(\Omega X)}\otimes\F{p}\cong UL$ for some
graded Lie algebra $L$~\cite{scott:00}, and indeed the proof of
Theorem~\ref{thm:bigyin} consists of a careful $p$-local lift
of this isomorphism.

Let $X \in \cw{n}{q}$.
Anick has constructed a differential Lie algebra $\Lie{X}$
and a chain algebra quasi-isomorphism 
$\theta_{X}:U\Lie{X}\rightarrow C_{*}(\Omega X;R)$
that commutes with the respective diagonals up to chain algebra
homotopy~\cite{anick:89}.
It follows that 
$\free{H(\theta_{X})}:\free{H(U\Lie{X})}\rightarrow 
	\free{H_{*}(\Omega X;R)}$
is an isomorphism of Hopf algebras.
Theorem~\ref{thm:bigyin} then follows immediately from the following
algebraic result, to be proved in Section~\ref{sec:free}.

\begin{THM}\label{thm:fhul}
Let $L$ be a connected, $R$-free differential Lie algebra
of finite type.
The inclusion $P \subset \free{H(UL)}$ of the sub-Lie algebra
of primitive elements of the Hopf algebra $\free{H(UL)}$ extends
to an isomorphism $UP \iso \free{H(UL)}$ of Hopf algebras.
\end{THM}

Our tool for spotting universal enveloping algebras is a
differential, torsion-free Andr\'e-Sj\"odin theorem.
In our context, a differential $\Gamma$-Hopf algebra is a differential
Hopf algebra which, as an algebra, is an algebra with divided powers.

\begin{THM}[c.f.~\cite{andre:71,sjodin:80}]\label{thm:andre-sjodin}
Let $(A,d)$ be a $R$-free, commutative differential Hopf algebra
of finite type.
Then $(A,d)^{\dual} \cong U(L,\partial)$ for some differential Lie algebra
$(L,\partial)$
if and only if 
$(A,d)$ is a differential $\Gamma$-Hopf algebra.
\end{THM}

The proof of Theorem~\ref{thm:andre-sjodin} is deferred to 
Section~\ref{sec:andre-sjodin}, but we will take the result as fact
throughout the rest of the paper.
 
Let $\varphi:\pi_{*}(\Omega X) \rightarrow H_{*}(\Omega X)$ be the
Hurewicz homomorphism.
A theorem of Cartan and Serre (see~\cite{f-h-t:01}) asserts that 
$\varphi\otimes\Q$ is
an isomorphism onto the subspace of primitives.  Along with the 
Milnor-Moore theorem, this gives the stunning result that
$U(\pi_{*}(\Omega X) \otimes \Q) \cong H_{*}(\Omega X;\Q)$.
The analogous torsion-free result over $R \subset \Q$ does not hold,
in general;  see Example~\ref{ex:K}. 
Instead we have the following theorem.

\begin{THM}\label{thm:hurewicz}
With the same hypotheses as in Theorem~\ref{thm:bigyin}, the induced
morphism
$\free{\varphi}:\free{\pi_{*}(\Omega X)\otimes R}
	\rightarrow \free{H_{*}(\Omega X;R)}$
maps injectively into $P$, and $P/\image{F\varphi}$ is torsion.
\end{THM}

A finite, simply-connected complex $X$ is called \emph{elliptic} if
$\pi_{*}(X)\otimes\Q$ is a finite-dimensional vector space:
$\sum_{m=1}^{\infty}\dim\pi_{m}(X)\otimes\Q < \infty$.
If $X$ is elliptic, then for sufficiently large primes $p$, 
$H_{*}(\Omega X_{(p)})$ is torsion-free and finitely generated as an
algebra~\cite{mcgibbon-wilkerson:86}.
Unfortunately, the notion of ``sufficiently large'' has not been made
precise.
As a corollary to Theorems~\ref{thm:bigyin} and~\ref{thm:hurewicz},
we obtain the following result along the same lines.

\begin{cor}
If $X \in \cw{n}{q}$ is elliptic, then $\free{H_{*}(\Omega X;R)}$ is 
finitely generated as an algebra, provided $\rho(R) \geq n/q$.
\end{cor}

\begin{proof}
For $m \geq 1$, 
$\dim(\pi_{m}(\Omega X)\otimes\Q) = \rank(\free{\pi_{m}(\Omega X)\otimes R})$.
As $X$ is rationally elliptic, the total rank of 
$\free{\pi_{*}(\Omega X)\otimes R}$
is finite.
By Theorem~\ref{thm:hurewicz}, the rank of 
$P \subset \free{H_{*}(\Omega X;R)}$
is also finite, and $P$ generates 
$\free{H_{*}(\Omega X;R)}$ by Theorem~\ref{thm:bigyin}.
\end{proof} 

The torsion module 
$P / \image{\free{\varphi}}$ is an obstruction to the existence
of a homotopy product decomposition of $\Omega X$,
localised at $p$, into the product of certain `atomic' spaces.
Let $\mathcal{C}$ be the collection of pointed spaces
$\{ S^{2m-1}, \Omega S^{2m+1} \}_{m \geq 1}
	\cup \{ U^{1} \} \cup \{ T^{2m+1}\{p^{r}\} \}_{m \geq 2, r \geq 1}
	\cup \{ S^{2m+1}\{p^{r}\} \}_{m \geq 1, r\geq 1}$.
Here $S^{2m+1}\{p^{r}\}$ is the homotopy-theoretic fibre of the degree
$p^{r}$ self-map on $S^{2m+1}$, and the spaces $U^{1}$,
$T^{2m+1}\{p^{r}\}$ are defined in~\cite{c-m-n:87}.
Denote by $\prod\mathcal{C}$ the collection of spaces having the weak
homotopy type of a finite-type product of spaces in $\mathcal{C}$.

\begin{thm}
If $X \in \cw{n}{q}$, $p \geq n/q$, and $\Omega X_{(p)} \in \prod\mathcal{C}$,
then $P/\image{\free{\varphi}}$ vanishes.
\end{thm}

\begin{proof}
Let $\varphi^{(r)}:E^{r}_{\pi}(\Omega X)\rightarrow E^{r}(\Omega X)$
be the morphism of Bockstein spectral sequences induced by the Hurewicz
morphism~\cite{neisendorfer:80}.
By~\cite[Theorem 9]{anick:92},
if $\Omega X_{(p)} \in \prod\mathcal{C}$, then
$E^{\infty}(\Omega X)$ is generated as an algebra by
$\image{\varphi^{(\infty)}}$.
But $E^{\infty}(\Omega X) \cong \free{H_{*}(\Omega X)}\otimes \F{p}
	\cong U(P \otimes \F{p})$,
so $P \otimes \F{p} \subset \varphi^{(\infty)}$.
We may identify $\varphi^{(\infty)} = \free{\varphi}\otimes\F{p}$,
so $p(P/\image{\free{\varphi}}) = P/\image{\free{\varphi}}$.
The result follows by Nakayama's Lemma.
\end{proof}
 
The author would like to thank Steve Halperin for suggesting the
problem, and Ran Levi for many helpful discussions.
Don Stanley spotted a small but critical error in an early draft of
the paper, and Peter Bubenik's careful proof-reading helped to make
many of the arguments presented here much clearer.

\section{Review}\label{sec:review}

\subsection{Notation and Conventions}

\emph{The ground ring $R$ will always be a principal 
ideal domain containing $1/2$}.
Algebraic objects are graded by the integers and are of finite type.
Algebras are augmented to the ground ring;  the augmentation ideal
is denoted $I(-)$.
Topological spaces are pointed and simply-connected.
  
The dual of a module $M$ is denoted by $M^{\dual}$.
Evaluation of $f \in M^{\dual}$ at $x \in M$
is denoted $\langle f,x \rangle$.
Set $\torsion{M} = \{ x \in M | rx = 0 \mathrm{\ for\ some\ } r \in R\}$.
The \emph{free part} of $M$ is the quotient module
$\free{M} = M / \torsion{M}$.
If $\varphi:M\rightarrow N$ is a linear map, then 
$\varphi(\torsion{M})\subset\torsion{N}$, so $\varphi$ factors to define
a linear map $\free{\varphi}:\free{M}\rightarrow\free{N}$.
Let $A$ be a differential Hopf algebra and let
$\kappa:H(A) \otimes H(A) \rightarrow H(A \otimes A)$ be the K\"unneth
homomorphism.
Then $\free{\kappa}$ is an isomorphism, and so
$\free{H(A)}$ is a Hopf algebra.

A \emph{Lie algebra}~\cite[Section 2]{c-m-n:79} is a non-negatively graded 
module $L = L_{\geq 0}$ along with a linear morphism
$[\, ,\,]:L \otimes L \rightarrow L$ called the Lie bracket,
that satisfies graded anti-symmetry, the graded Jacobi identity,
and $[x,[x,x]]=0$ if $\deg{x}$ is odd.
Since $\torsion{L}$ is a Lie ideal,
$\free{L}$ is a Lie algebra.
In fact, if $R = \Zloc{3}$, and $L$ is a torsion-free
module with a bilinear bracket satisfying anti-symmetry and the Jacobi
identity, then the last
condition is automatically satisfied
and so $L$ is a Lie algebra.
In particular if $X$ is a space, then 
$\free{\pi_{*}(\Omega X)\otimes\Zloc{3}}$ is a 
Lie algebra even though $\pi_{*}(\Omega X)\otimes\Zloc{3}$ may not be
(the canonical example is $\pi_{*}(\Omega S^{2n})$~\cite{liulevicius:62}).
If $A$ is an algebra, then the commutator bracket
$[a,a'] = aa' - (-1)^{\deg{a}\deg{a'}}a'a$ makes $A$ into a Lie algebra.
A Lie algebra is called \emph{connected} if it is concentrated in strictly
positive degrees.
A \emph{differential Lie algebra} is a chain complex $(L,\partial)$,
where $L$ is a Lie algebra and the differential $\partial$ satisfies
the Leibniz condition
$\partial[x,y] = [\partial{x}, y] + (-1)^{\deg{x}}[x,\partial{y}]$.
The homology of a differential Lie algebra is a Lie algebra.
The \emph{universal enveloping algebra} of a Lie algebra
$L$ is an associative algebra $UL$ along with a Lie algebra morphism
$\iota:L \rightarrow UL$ such that, if
$A$ is an algebra and $f:L \rightarrow A$ is a Lie algebra morphism,
then there exists a unique algebra morphism
$F:UL \rightarrow A$ such that $F \circ \iota = f$.
It follows that the universal enveloping algebra of a Lie algebra
is unique up to canonical algebra isomorphism.
The Lie algebra morphism $L \rightarrow UL \otimes UL$
defined by $x \mapsto x \otimes 1 + 1 \otimes x$ extends to define a
diagonal $\Delta:UL \rightarrow UL \otimes UL$
that provides $UL$ with a natural cocommutative Hopf algebra structure.
If $L$ is a differential Lie algebra, then the differential extends
to a derivation on $UL$ making it a differential Hopf algebra.

A $\Gamma$-\emph{algebra} is a connected, non-negatively graded commutative 
algebra $A=A^{\geq 0}$, along with a system of \emph{divided powers}
$\gamma^{k}:A^{2n} \rightarrow A^{2nk}$ for $k \geq 0$, $n \geq 1$, that
satisfy a list of axioms~\cite{andre:71}.
We will need explicitly the following axiom:
\begin{eqnarray*}
	\gamma^{j}(a)\gamma^{k}(a) = \binom{j+k}{j}\gamma^{j+k}(a);
	&	a \in A^{2n},\; j, k, n \geq 1 .
\end{eqnarray*}
By induction it follows that $a^{k} = k!\gamma^{k}(a)$ for $a \in A^{2n}$
and $k, n \geq 1$.
Henceforth we assume that all $\Gamma$-algebras are augmented.
A $\Gamma$-\emph{morphism} is an algebra morphism $\varphi:A \rightarrow B$,
where $A$ and $B$ are $\Gamma$-algebras and
$\varphi(\gamma^{k}(a))=\gamma^{k}(\varphi(a))$ for all
$a \in IA^{\mathrm{even}}$ and $k \geq 0$.
A $\Gamma$-\emph{derivation} $\theta$ on a $\Gamma$-algebra $A$ is an algebra
derivation that satisfies the additional condition
$\theta\gamma^{k}(a) = \theta(a)\gamma^{k-1}(a)$ for 
$a \in IA^{\mathrm{even}}$.
A \emph{differential $\Gamma$-algebra} is a cochain algebra
$(A,d)$, where $A$ is a $\Gamma$-algebra and $d$ is a $\Gamma$-derivation.
A $\Gamma$-\emph{Hopf algebra} is simultaneously a $\Gamma$-algebra
and a Hopf algebra, such that the diagonal map is a $\Gamma$-morphism.
A \emph{differential $\Gamma$-Hopf algebra} is a differential Hopf algebra 
$(A,d)$ where $A$ is a $\Gamma$-Hopf algebra and $d$ is a $\Gamma$-derivation.

\subsection{The Bockstein spectral sequence}\label{sec:bss}

We review the relevant details of the Bockstein
spectral sequence associated to the short exact coefficient sequence 
$0 \rightarrow \Zloc{p} \xrightarrow{\times p} \Zloc{p}
	\rightarrow \F{p} \rightarrow 0$.
The standard reference is Browder~\cite{browder:61};  for the
particulars of the Bockstein spectral sequence of the universal enveloping
algebra of a differential Lie algebra, see~\cite{scott:00}.

For a $\Zloc{p}$-free chain complex $C$, the Bockstein spectral sequence
is a spectral sequence of chain complexes over $\F{p}$.
If $H(C)$ is finite type, then the spectral sequence converges,
$H(C;\F{p}) \Rightarrow \free{H(C)}\otimes \F{p}$.

Let $(L,\partial)$ be a connected, $\Zloc{p}$-free differential Lie 
algebra.
Then by Theorem~\ref{thm:andre-sjodin}, $(UL)^{\dual}$ is a differential 
$\Gamma$-Hopf algebra.
For brevity, set $G = (UL)^{\dual}$, and let $\{(E_{r},\beta_{r})\}$ be the 
Bockstein spectral sequence modulo $p$ for $G$.
For $r \geq 1$, let $\rho_{r}:H(G) \rightarrow E_{r}$ be the usual
morphism~\cite[Section 1]{browder:61}.  
Each $\rho_{r}$ is an algebra morphism.
For each $r \geq 1$ there is an injective differential $\Gamma$-morphism 
$m_{r+1}:(E_{r+1},0) \qi (E_{r},\beta_{r})$ that induces an isomorphism
in homology, and is used to identify $E_{r+1} \iso H(E_{r})$ as a 
$\Gamma$-Hopf algebra.
It follows that for $\zeta \in H(G)$,
\begin{equation}\label{eq:bss-maps}
	m_{r+1}\circ\rho_{r+1}(\zeta) = \rho_{r}(\zeta) + \beta_{r}(y_{r})
\end{equation}
for some $y_{r} \in E_{r}$.
For $1 \leq r \leq s$, define injective $\Gamma$-morphisms
$\sigma_{rs}:E_{s} \rightarrow E_{r}$ 
by $\sigma_{rr} = 1_{r}:E_{r}\rightarrow E_{r}$ for $r \geq 1$ and
$\sigma_{rs} = m_{s}\circ\cdots\circ m_{r+1}$ for $1 \leq r < s$.
The sequence $\{ E_{r} \}$ then forms an inverse system,
and $E_{\infty} = \varprojlim{E_{r}}$.
Using this characterisation, one can show that
$E_{\infty}$ is a $\Gamma$-Hopf algebra;  in particular,
by~\cite{sjodin:80}, $E_{\infty} \iso \Gamma(V_{\infty})$ as a
$\Gamma$-algebra.
The map $\sigma_{\infty} = \varprojlim{\sigma_{1r}}$ is an injective 
$\Gamma$-morphism.
Finally, $\rho_{\infty} = \varprojlim{\rho_{r}}:H(G) \rightarrow E_{\infty}$
is a surjective algebra morphism, and
$\ker{\rho} = \torsion(H(G)) + pH(G)$.

\section{Homology of $UL$ modulo torsion}\label{sec:free}

In this section we prove the main algebraic result of the paper,
Theorem~\ref{thm:fhul} from the introduction, and discuss some properties of 
the morphism $\free{H(L)}\rightarrow\free{H(UL)}$ induced by the canonical 
inclusion $L \rightarrow UL$.
Since $UL$ is a differential Hopf algebra, $FH(UL)$ is a Hopf algebra.

\begin{proof}[Proof of Theorem~\ref{thm:fhul}]
Let $P \subset FH(UL)$ be the submodule of primitive elements. 
The inclusion extends to a Hopf algebra morphism
$\beta:UP \rightarrow FH(UL)$.
Recall that $\beta$ is an isomorphism if and only if
$\beta_{\mathfrak{p}}:(UP)_{\mathfrak{p}}\rightarrow FH(UL)_{\mathfrak{p}}$
is an isomorphism for all prime ideals $\mathfrak{p}$ in 
$R$~\cite{atiyah-macdonald}.
Furthermore, each prime ideal $\mathfrak{p}$ in $R$ is of the form
$pR$ for some prime $p \geq \rho(R)$, and
$R_{\mathfrak{p}} \cong \Zloc{p}$.
Therefore, without loss of generality, we may
assume that $L$ is a differential Lie algebra over $\Zloc{p}$,
$p \geq \rho(R)$.

To begin, we construct a subalgebra $A$ of $H[(UL)^{\dual}]$
such that $A \otimes \F{p} \iso E_{\infty}$.
Recall (Section~\ref{sec:bss}) that $E_{\infty} \iso \Gamma(V_{\infty})$.
Let $\{ v_{i} \}$ be a basis of $V_{\infty}$.
Since $\rho_{\infty}$ is surjective, each $v_{i} = \rho_{\infty}(\zeta_{i})$ 
for some $\zeta_{i} \in H(G)$.
Let $n = \deg{v_{i}}$.
There exists $N > 0$ such that
$\beta_{r}(E^{\leq n}_{r}) = 0$ whenever $r \geq N$.
It follows from (\ref{eq:bss-maps}) that
$\sigma_{\infty}(v_{i}) = \rho_{1}(\zeta_{i} + \sum_{r=1}^{N}\xi_{r})$
where $\xi_{r} \in H(G)$ has order $p^{r}$.
Let $z_{i} \in G$ be a cycle whose homology class is
$\zeta_{i} + \sum_{r=1}^{N}\xi_{r}$.
Let $V$ be the free $\Zloc{p}$-module on the basis $\{ v_{i} \}$.
Define a linear map $\sigma: V \rightarrow G$ by
$\sigma(v_{i}) = z_{i}$.
Then $\sigma$ extends to a differential $\Gamma$-morphism
$A = (\Gamma(V),0) \xrightarrow{\sigma} G$, and 
$\sigma \otimes \F{p} = \sigma_{\infty}$.
Since $\sigma_{\infty}$ is injective and $\Gamma(V)$, $G$ are
$\Zloc{p}$-free as modules, $\sigma$ is injective.

Let $\pi : A \rightarrow E_{\infty}$ be reduction mod $p$.
Passing to homology, it is a routine exercise to verify that
$\rho_{1}\circ\sigma_{*} = \sigma_{\infty}\circ\pi$ and
$\rho_{\infty}\circ\sigma_{*} = \pi$.
Suppose $\sigma_{*}(a) = 0$ for some $a \in A$.
Then $\pi(a) = \rho_{\infty}\circ\sigma_{*}(a) = 0$, so
$a = pa'$ for some $a' \in A$.
Thus $p \sigma_{*}(a') = 0$, so $\sigma_{*}(a')$ is a torsion
element of $H(G)$.
It follows that $\rho_{\infty}\circ\sigma_{*}(a') = 0$,
whence $\pi(a') = 0$.
So $a' = pa''$ for some $a'' \in A$, which is again
sent to a torsion element of $H(G)$ by $\sigma_{*}$.
Repeat the argument;  eventually the process ends since
$A$ has no infinitely $p$-divisible elements.
Therefore $a=0$ and so $\sigma_{*}$ is injective.

Let $T$ be the torsion ideal of $H(G)$.
Since $T \subset \ker{\rho_{\infty}}$, the morphism $\rho_{\infty}$ factors
to define $\bar{\rho}_{\infty}:\free{H(G)} \rightarrow E_{\infty}$.
Let $\hat{\sigma}_{*}$ be the composite
$A \xrightarrow{\sigma_{*}} H(G) \rightarrow \free{H(G)}$.
One verifies that $\bar{\rho}_{\infty}\circ\hat{\sigma}_{*} = \pi$.
From the definitions, $\hat{\sigma}_{*} \otimes \F{p}$
is the identity on $E_{\infty}$.
Thus 
$\free{H(G)}/\image{\hat{\sigma}_{*}}=p(\free{H(G)}/\image{\hat{\sigma}_{*}})$.
Since $\free{H(G)}/\image{\hat{\sigma}_{*}}$ is degree-wise finitely
generated,
Nakayama's Lemma states that $\free{H(G)}/\image{\hat{\sigma}_{*}}=0$.
Therefore
$\hat{\sigma}_{*}$ is surjective and hence an isomorphism of algebras.

By Theorem~\ref{thm:andre-sjodin}, $G$ is a differential $\Gamma$-Hopf
algebra.
We claim that the diagonal $\psi$ in $G$ 
induces a $\Gamma$-Hopf algebra structure on $A$.
Indeed, since $\sigma$ is a $\Gamma$-morphism, 
$\sigma_{*}(\gamma^{k}(a)) = [\gamma^{k}(\sigma(a))]$.
A $\Gamma$-morphism
$\chi:A \rightarrow A \otimes A$ is defined by
$(\hat{\sigma}_{*}^{-1} \otimes \hat{\sigma}_{*}^{-1})\circ
	\free{H(\psi)}\circ\hat{\sigma}_{*}$.
Since $\psi$ is a $\Gamma$-morphism, so too is $\chi$.
Furthermore, $\chi$ inherits coassociativity and the counit from
$\psi$, and so $\chi$ provides $A$ with a $\Gamma$-Hopf algebra structure.

By definition, $\chi$ commutes with $\sigma$ ``up to torsion''.
That is, given $a \in A$, there exist 
$\Phi, \Psi \in G$, and $r \geq 0$ such that
\begin{equation}\label{eq:chi}
	(\sigma \otimes \sigma)\circ\chi(a) = \psi\circ\sigma(a) + \Phi
\end{equation}
and $d\Psi = p^{r}\Phi$.

We now define a pairing
\begin{equation}\label{eq:pair}
	A \otimes \free{H(UL)} \xrightarrow{\langle\,,\,\rangle} \Zloc{p}
\end{equation}
that exhibits $\free{H(UL)}$ 
as the Hopf algebra dual of the $\Gamma$-Hopf algebra $A$.
By Theorem~\ref{thm:andre-sjodin}, the proof will then be complete.
Use $\sigma$ to identify $(A,0)$ as a differential $\Gamma$-subalgebra
of $G$ to obtain a pairing
$A \otimes \ker{\partial} \xrightarrow{\langle\,,\,\rangle} \Zloc{p}$,
where $\partial$ is the differential on $UL$.
Let $T = \{ x \in UL | p^{r}x \in\image\partial
	\mathrm{\ for\ some\ }r \geq 0\}$. 
Since the differential in $A$ vanishes, $\langle A,\image\partial \rangle = 0$.
Furthermore, $\Zloc{p}$ is a domain, hence $\langle A, T \rangle = 0$.
Therefore $\langle\,,\,\rangle$ factors to define~(\ref{eq:pair}).

Observe that (\ref{eq:pair}) is non-singular.
Indeed, as a chain map, $\sigma:(A,0) \rightarrow G$ splits, and $\sigma$ 
induces an isomorphism $A \iso \free{H(G)}$. 
It follows that $G = A \oplus B \oplus C$, where $d : B \rightarrow C$ and
$H(G) = A \oplus (C/dB)$.
Dually, $UL = A^{\dual} \oplus B^{\dual} \oplus C^{\dual}$, and
$H(UL) = A^{\dual} \oplus (B^{\dual}/\partial C^{\dual})$, where
$B^{\dual}/\partial C^{\dual}$ is the torsion submodule.
Therefore $\free{H(UL)} \iso A^{\dual}$ as a module, and one can check
that the duality is provided by (\ref{eq:pair}).

Since $\sigma$ is an algebra morphism,
$\langle a \otimes a',\Delta x \rangle = \langle aa',x \rangle$
for $x \in UL$ and $a,a' \in A$.
Pass to (\ref{eq:pair}) to see that the multiplication
in $A$ is dual to the diagonal $\free{H(\Delta)}$.

Suppose that $x,y \in UL$ are cycles, and $a \in A$.
Using~(\ref{eq:chi}), we see that
$\langle a, xy \rangle = \langle \psi(a), x \otimes y \rangle
	= \langle \chi(a), x \otimes y \rangle + \langle \Phi, x\otimes y \rangle$,
where $p^{r}\Phi = d\Psi$ for some $\Psi \in UL \otimes UL$.
Since $p^{r}\Phi$ is a boundary, $x$ and $y$ are cycles, and $\Zloc{p}$
is an integral domain, it follows that
$\Phi(x\otimes y) = 0$.
Pass to~(\ref{eq:pair}) to conclude that $\chi$ is dual to the
multiplication in $\free{H(UL)}$.
\end{proof}

\begin{prop}\label{prop:Lie-embed}
The natural map $\iota:L\rightarrow UL$ induces an injection of
Lie algebras $\free{\iota_{*}}:\free{H(L)} \hookrightarrow P$,
with $P / \image{\free{\iota_{*}}}$ torsion.
\end{prop}

\begin{proof}
Let $\bar{\Delta}$ be the reduced diagonal on $UL$,
and let $\kappa:H(UL)\otimes H(UL) \rightarrow H(UL\otimes UL)$ be
the K\"unneth morphism.
Recall that $\free{\kappa}$ is an isomorphism.
Since $\bar{\Delta}\circ\iota = 0$, 
$(\free{\kappa})^{-1}\circ\free{\bar{\Delta}_{*}}\circ\free{\iota}_{*} = 0$,
and so $\image \free{\iota}_{*} \subset P$.
Observe that $\free{H(L)}\otimes\Q = H(L\otimes\Q)$
and $\free{H(UL)}\otimes\Q = H(U(L\otimes\Q)) = UH(L\otimes\Q)$.
Furthermore, $\free{\iota_{*}}\otimes\Q$ is the canonical
inclusion $H(L\otimes\Q)\rightarrow UH(L\otimes\Q)$.
As $\free{H(L)}$ and $\free{H(UL)}$ are torsion-free, the result now follows.
\end{proof}

\begin{ex}
Define a differential Lie algebra $L$ as follows.
As a chain complex, $L$ is the free $\Zloc{3}$-module on the graded basis
$\{ x_{1} , y_{2} , x_{3}, x_{5} \}$,
with subscript indicating degree;
the only non-zero differential is $dy_{2} = 3 x_{1}$.
The bracket is defined by $[ y_{2} , x_{1} ] = x_{3}$,
$[ y_{2} , x_{3} ] = x_{5}$, and all others vanishing.
Let $\langle S \rangle$ denote the abelian Lie algebra on the graded set $S$.
Then
$\free{H(L)}= \langle [x_{3}], [x_{5}] \rangle$.

By computing the Bockstein spectral sequence for $UL$, one finds that
$\free{H(UL)} = U \langle z_{3}, z_{5} \rangle$,
where $z_{3}$ and $z_{5}$ are represented in $UL$ by
$x_{3}$ and $x_{1} y_{2}^{2} + y_{2}x_{3}$, respectively.
Let $\iota:L\rightarrow UL$ be the canonical injection.
Clearly $\free{\iota_{*}}[x_{3}]=z_{3}$;
however,
$dy_{2}^{3} = -6 x_{5} + 9(x_{1}y_{2}^{2} + y_{2}x_{3})$ in $UL$,
and so $\free{\iota_{*}}[x_{5}] = (3/2)z_{5}$.
Therefore, $P / \image\free{\iota_{*}} \iso \Z/3\Z$, with
generator represented by $z_{5}$.
\end{ex}

\section{The Hurewicz homomorphism}\label{sec:hurewicz}

In this section we prove Theorem~\ref{thm:hurewicz}.
Suppose $X \in \cw{n}{q}$.
As usual, set $P = P(FH_{*}(\Omega X;R))$.

We will need the following lemma.
Let $Q \subset H_{*}(\Omega X;\Q)$ be the sub Lie algebra of primitive
elements.
The inclusion extends to an isomorphism of Hopf algebras
$\beta:UQ \iso H_{*}(\Omega X;\Q)$~\cite{milnor-moore:65}.

\begin{lem}\label{lem:P-and-Q}
$P \otimes \Q \iso Q$.
\end{lem}

\begin{proof}
By universal coefficients, the natural Hopf algebra morphism
$\mu:H_{*}(\Omega X;R)\otimes\Q \rightarrow H_{*}(\Omega X;\Q)$, defined by
$\mu([z] \otimes x) = [z \otimes x]$ for $x \in \Q$ and cycles
$z \in C_{*}(\Omega X;R)$, is an isomorphism.
Let $r:H_{*}(\Omega X;R) \rightarrow \free{H_{*}(\Omega X;R)}$ be the 
quotient morphism.
Since 
$r\otimes\Q:H_{*}(\Omega X;R)\otimes\Q
	\rightarrow[\free{H_{*}(\Omega X;R)}]\otimes\Q$
is an isomorphism, we get an isomorphism of Hopf algebras
$\bar{\mu}:[\free{H_{*}(\Omega X;R)}]\otimes\Q
	\xrightarrow{\iso} H_{*}(\Omega X;\Q)$
defined by $\bar{\mu}=\mu\circ(r\otimes\Q)^{-1}$.
By Theorem~\ref{thm:bigyin} and the mapping property of universal enveloping 
algebras, 
$\alpha\otimes\Q:U(P\otimes\Q)\xrightarrow{\iso}
	[\free{H_{*}(\Omega X;R)}]\otimes\Q$.
It follows that $\beta^{-1}\circ\bar{\mu}\circ(\alpha\otimes\Q)$ 
restricts to an
isomorphism of Lie algebras $P \otimes \Q \xrightarrow{\iso} Q$.
\end{proof}

\begin{proof}[Proof of Theorem~\ref{thm:hurewicz}]
Let
$\varphi:\pi_{*}(\Omega X) \rightarrow H_{*}(\Omega X)$
be the Hurewicz homomorphism.
The inclusion $R \subset \Q$ induces the commutative
diagram
\[
\xymatrix{
	\pi_{*}(X)\otimes R \ar[r]^{\varphi\otimes R} \ar[d]
		& H_{*}(\Omega X;R) \ar[d]	\\
	\pi_{*}(\Omega X)\otimes\Q \ar[r]_{\varphi\otimes\Q}
		& H_{*}(\Omega X ; \Q)
}
\]
in which the right arrow is a Hopf algebra morphism, and the other
arrows preserve brackets.
The kernels of the vertical arrows are precisely the respective torsion
submodules.
Therefore the diagram
\[
\xymatrix{
	\free{\pi_{*}(X)\otimes R} \ar[r]^{\free{\varphi\otimes R}} \ar[d]	
		& \free{H_{*}(\Omega X;R)} \ar[d] 
		& UP \ar[l]_(0.4){\iso}^(0.4){\alpha} \ar[d]	\\
	\pi_{*}(\Omega X)\otimes\Q \ar[r]_{\varphi\otimes\Q}
		& H_{*}(\Omega X ; \Q) 
		& U(P \otimes \Q) \ar[l]_{\iso}^{\beta}
}
\]
commutes.
The left vertical arrow, $\free{\varphi\otimes R}$, and $\varphi\otimes\Q$ are
Lie algebra morphisms.
The remaining arrows are Hopf algebra morphisms.
The vertical arrows are injections,
$\varphi\otimes\Q$ is an isomorphism onto $P \otimes \Q$ by 
Lemma~\ref{lem:P-and-Q}, and the right vertical arrow
is of the form $U\{ P \rightarrow P \otimes \Q \}$.
It follows that 
$\free{\varphi\otimes R}:\free{\pi_{*}(X)\otimes R}\rightarrow P$ is an
injection and that $P / \image\free{\varphi\otimes R}$ is torsion.
\end{proof}

\begin{ex}\label{ex:K}
For this example, all spaces are localised at a fixed prime $p$.
Let $\alpha:S^{2p}\rightarrow S^{3}$ be a representative of the generator
of $\pi_{2p}(S^{3}) \cong \Z/p\Z$.
Set
$K = S^{4} \cup_{\Sigma\alpha} e^{2p+2}$.
The map $\Sigma\alpha\circ p:S^{2p+1}\rightarrow S^{4}$ is null-homotopic, 
and so we get a map $f:S^{2p+2}\rightarrow K$ making 
the diagram
\[
\xymatrix{
	S^{2p+1} \ar[r] \ar[d]_{p}	& {*} \ar[d] \ar[r]
		& S^{2p+2} \ar[d]^{f} \\
	S^{2p+1} \ar[r]^{\Sigma\alpha}	& S^{4} \ar[r]	& K
}
\]
commute up to homotopy, where the rows are cofibration
sequences.
From the associated long exact sequence in homology, we deduce that
in degree $2p+2$, $f_{*}$ is multiplication by $p$.

By the Bott-Samelson theorem, 
$H_{*}(\Omega K) \cong T( e_{3}, e_{2p+1} )$ and
$P = PH_{*}(\Omega K) = L( e_{3}, e_{2p+1} )$.
Let $\varphi:\pi_{*}(\Omega K) \rightarrow H_{*}(\Omega K)$ be the 
Hurewicz map.
Then $\varphi(f^{\sharp}) = p e_{2p+1}$, where
$f^{\sharp}:S^{2p+1}\rightarrow \Omega K$ is the adjoint of $f$.
Since $\varphi(f^{\sharp}) \neq 0$ and $H_{*}(\Omega K)$ is torsion-free,
$f^{\sharp}$ is not a torsion element.
It follows that $p e_{2p+1} \in \image{\free{\varphi}}$.

If $e_{2p+1} \in \image{\free{\varphi}}$, then 
$e_{2p+1} \in \image{\varphi}$.
By Baues~\cite[Lemma V.3.10]{baues:81}, $\Omega K$ would be homotopy-equivalent
to the product of odd-dimensional spheres and loops on the same.
But $K$ is a retract of $\Sigma^{2}\mathbf{C}P^{p}$, and so
$H^{*}(\Omega K;\F{p})$ has a non-vanishing $\mathcal{P}^{1}$.
Therefore $e_{2p+1} \not\in \image{\free{\varphi}}$, and so
$P / \image{\free{\varphi}} \cong \Z/p\Z$.
\end{ex}

\section{A torsion-free Andr\'e-Sj\"odin theorem}\label{sec:andre-sjodin}

This section is devoted to the proof of Theorem~\ref{thm:andre-sjodin}.
Let $A$ be an $R$-free $\Gamma$-algebra.
The module $\Hom(A,A)$ is an algebra under composition, and hence
is a Lie algebra with the commutator bracket.
The submodule $\Der{A}$ of all
$\Gamma$-derivations on $A$ is a sub-Lie algebra of $\Hom(A,A)$.

We will need a construction of Gulliksen and 
Levin~\cite{gulliksen-levin:69}.
Suppose $A$ is an $R$-free $\Gamma$-Hopf algebra,
with diagonal $\Delta$.
Given $f \in A^{\dual}$, define $\nu(f)$ to be the composite
\[
	A \xrightarrow{\Delta} A \otimes A \xrightarrow{1 \otimes f}
		A \otimes R = A.
\]
Define
$DA = IA \cdot IA + \sum_{k \geq 2} \gamma^{k}(IA^{\mathrm{even}})$.
Set $\tilde{P}(A^{\dual}) = \{ f \in A^{\dual} | f(DA) = 0 \}$.
The primitives of $A^{\dual}$ may be identified as
$P(A^{\dual}) = 
	\{ f \in A^{\dual} | f(IA\cdot IA) = 0 \}$~\cite{milnor-moore:65}.

\begin{lem}\label{lem:sjodin}\cite{sjodin:80}
With the notation and hypotheses above,

(a) $\nu:A^{\dual}\rightarrow \Hom(A,A)$ is an algebra morphism;

(b) $\nu(f)|_{A^{\deg{f}}} = f$;

(c) $\nu(f) \in \Der{A}$ if and only if $f \in \tilde{P}(A^{\dual})$;

(d) $\tilde{P}(A^{\dual})$ is a sub-Lie algebra of $A^{\dual}$ and
$\nu: \tilde{P}(A^{\dual}) \rightarrow \Der{A}$ is a Lie monomorphism.
\end{lem}

\begin{prop}\label{prop:g-hopf}
Let $A$ be a $\Zloc{p}$-free $\Gamma$-Hopf algebra over $\Zloc{p}$.

(a) There is a submodule $X \subset A$ whose inclusion extends to an
isomorphism of $\Gamma$-algebras
$\Gamma(X) \xrightarrow{\cong} A$.

(b) The inclusion $P(A^{\dual}) \subset A^{\dual}$ extends to an isomorphism
of Hopf algebras $UP(A^{\dual}) \xrightarrow{\cong} A^{\dual}$. 
\end{prop}

\begin{proof}
(a) Set $\bar{A} = A \otimes \F{p}$.
Then $\bar{A}$ is a $\Gamma$-Hopf algebra over $\F{p}$.
Let $\sigma:\bar{A}/D\bar{A} \rightarrow \bar{A}$ be a splitting of
the quotient map.
By Sj\"odin~\cite{sjodin:80}, $\sigma$ extends to an isomorphism of
$\Gamma$-algebras
$\Gamma(\bar{A}/D\bar{A}) \xrightarrow{\cong} \bar{A}$.
Let $\{ \bar{x}_{i} \}$ be a basis of $\image{\sigma} \subset \bar{A}$
indexed over the positive integers.
Choose a set of representatives $\{ x_{i} \}$ in $A$, and let
$X$ be the linear span.
The inclusion $X \subset A$ extends to a $\Gamma$-morphism
$\varphi:\Gamma(X) \rightarrow A$.
The reduction mod $p$, 
$\varphi\otimes\F{p}:\Gamma(\image{\sigma})\rightarrow\bar{A}$,
is an isomorphism.
In particular, $\varphi\otimes\F{p}$ is surjective, so by Nakayama's Lemma,
$A/\image{\varphi}=0$.
It follows that $\varphi$ is surjective.
Since $A$ and $\Gamma(X)$ are $\Zloc{p}$-free and finite type, and 
$\varphi\otimes\F{p}$ is injective, $\varphi$ is itself injective
and therefore an isomorphism.

(b) To begin, we claim that $P(A^{\dual}) = \tilde{P}(A^{\dual})$.
From the definitions, $\tilde{P}(A^{\dual}) \subset P(A^{\dual})$.
It suffices to show that if $f \in P(A^{\dual})$, then
$f(\gamma^{k}(a))=0$ for all $k \geq 2$ and $a \in IA^{\mathrm{even}}$.
Since $a^{k} = k!\gamma^{k}(a)$ and $f(IA \cdot IA) = 0$,
$0 = f(a^{k}) = k!f(\gamma^{k}(a))$.
As $\Zloc{p}$ is a domain, the claim follows.

To show that the canonical algebra morphism 
$UP(A^{\dual})\rightarrow A^{\dual}$ is surjective, we use an argument
of Sj\"odin.
As a module, $IA \cong X \oplus DA$.
For $i \geq 1$, define $f_{i} \in A^{\dual}$ by
$f_{i}(\Zloc{p} \oplus DA) = 0$ and $f_{i}(x_{j}) = \delta_{ij}$.
By definition, each $f_{i} \in P(A^{\dual})$.
Let $\mathcal{A}$ be the set of finite sequences of non-negative integers
$r = (r_{i})$ such that $r_{i}=0$ or $1$ if $\deg{x_{i}}$ is odd,
and the last term in $r$ is non-zero.
Define $n(r)$ to be the length of $r$ (so $r_{n(r)} \neq 0$).
Order $\mathcal{A}$ by setting $r > r'$ if $n(r) > n(r')$ or
$n(r)=n(r')$ and $r_{n(r)}>r_{n(r')}$.
Set $f^{r} = f_{1}^{r_{1}}\cdots f_{n}^{r_{n}}$ and
$\gamma^{r}(x) = \gamma^{r_{n}}(x_{n})\cdots\gamma^{r_{1}}(x_{1})$
where $n = n(r)$.
Since $\nu$ is an algebra morphism,
$\nu(f^{r}) = \nu(f_{1})^{r_{1}}\circ\cdots\circ\nu(f_{n})^{r_{n}}$
(exponents with respect to composition).
As $f_{i} \in P(A^{\dual}) = \tilde{P}(A^{\dual})$,
$\nu(f_{i})$ is a $\Gamma$-derivation;  that is,
$\nu(f_{i})(\gamma^{k}(x_{i})) = \nu(f_{i})(x_{i})\gamma^{k-1}(x_{i})$.
It follows that $\nu(f^{r})(\gamma^{r}(x))=1$,
and $\nu(f^{r})(\gamma^{r'}(x))=0$ if $r>r'$.
By Lemma~\ref{lem:sjodin}(b), $f^{r}(\gamma^{r}(x))=1$, while
$f^{r}(\gamma^{r'}(x))=0$ if $r>r'$ and $\deg{f^{r}}=\deg{\gamma^{r'}(x)}$.
If $r>r'$ and $\deg{f^{r}}\neq\deg{\gamma^{r'}(x)}$, then of course
$f^{r}(\gamma^{r'}(x))=0$ for degree reasons.
The set $\{\gamma^{r}(x)\}_{r\in\mathcal{A}}$ forms an additive basis for
$A$.
The above formulae for the $f^{r}$ imply that $\{ f^{r} \}_{r\in\mathcal{A}}$
is a basis for $A^{\dual}$.
In particular, the set $\{ f_{i} \}$ generates $A^{\dual}$ as an algebra,
and so $UP(A^{\dual})\rightarrow A^{\dual}$ is surjective.

Let $\chi_{V}$ denote the Euler-Poincar\'e series of the non-negatively
graded $\F{p}$-space $V$.
We extend the definition to a non-negatively graded,
free $\Zloc{p}$-module $M$:
$\chi_{M} = \sum_{n\geq 0} (\rank{M_{n}})t^{n}$.
Note that $\chi_{M} = \chi_{M \otimes \F{p}}$.
Since $UP(A^{\dual})\rightarrow A^{\dual}$ is surjective,
$\chi_{UP(A^{\dual})} \geq \chi_{A^{\dual}}$.
If $f \in P(A^{\dual})$, then $(f\otimes\F{p})(D\bar{A})=0$,
so $P(A^{\dual})\otimes\F{p}\subset \tilde{P}(\bar{A}^{\dual})$
and the natural morphism 
$U(P(A^{\dual})\otimes\F{p})\rightarrow U\tilde{P}(\bar{A}^{\dual})
	=\bar{A}^{\dual}$
is an injection.
It follows that $\chi_{UP(A^{\dual})} \leq \chi_{A^{\dual}}$,
so in fact the Euler-Poincar\'e series are equal.
Therefore $UP(A^{\dual})\rightarrow A^{\dual}$ is an isomorphism.
\end{proof}

Before finally proving Theorem~\ref{thm:andre-sjodin}, we state a Lemma.

\begin{lem}\label{lem:g-alg}
Let $R$ be a characteristic-zero principal ideal domain.
Let $A$ and $B$ be $R$-free $\Gamma$-algebras.

(1) If $\varphi:A\rightarrow B$ is an algebra morphism, then it is a
$\Gamma$-morphism.

(2) If $\theta:A \rightarrow B$ is a derivation, then it is a
$\Gamma$-derivation.

In particular, if $A$ is a $\Gamma$-algebra and a differential Hopf algebra,
then $A$ is a differential $\Gamma$-Hopf algebra.
\end{lem}

\begin{proof}
Straightforward.
\end{proof}

\begin{proof}[Proof of Theorem~\ref{thm:andre-sjodin}]
Suppose $A$ is a differential $\Gamma$-Hopf algebra.
As usual, denote by $P(A^{\dual})$ the sub-Lie algebra of primitive
elements.
The natural Hopf algebra morphism $UP(A^{\dual})\rightarrow A^{\dual}$
is an isomorphism if and only if 
$UP(A^{\dual})_{\mathfrak{p}}\rightarrow A^{\dual}_{\mathfrak{p}}$
is an isomorphism for every prime ideal $\mathfrak{p}$ of $R$,
and that each prime ideal $\mathfrak{p}$ is of the form
$pR$ for some prime $p \geq \rho(R)$.
Furthermore, observe that $P(A^{\dual})_{(p)}=P(A^{\dual}_{(p)}$.
Therefore, without loss of generality we may assume that
$R = \Zloc{p}$ for some odd prime $p$. 
By Proposition~\ref{prop:g-hopf}, 
$A^{\dual} \cong UP(A^{\dual})$ as a Hopf algebra.
Since the differential in $A^{\dual}$ is a Hopf derivation,
it preserves $P(A^{\dual})$.
Thus $P(A^{\dual})$ forms a differential Lie algebra, establishing the
`if' assertion.

Conversely, suppose $(A,d)^{\dual} \cong U(L,\partial)$.
By Lemma~\ref{lem:g-alg}, it suffices to show that $A^{\dual}$ is a
$\Gamma$-algebra.
In~\cite{halperin:92}, Halperin constructed a chain isomorphism
$\gamma_{L}:UL \rightarrow (\Gamma(L^{\dual}),\bar{D})^{\dual}$,
where $\bar{D}$ is a $\Gamma$-derivation.
By construction, $\gamma_{L}$ factors as
\[
	UL \xrightarrow{U\sigma} U\mathcal{D} \xrightarrow{\theta}
		(\Gamma(L^{\dual}),\bar{D})^{\dual}
\]
for some differential Lie morphism $\sigma:L\rightarrow\mathcal{D}$.
By the proof of~\cite[Proposition 4.2]{halperin:92},
$\theta$ is a coalgebra morphism.
Since $U\sigma$ is also a coalgebra morphism, $\gamma_{L}$ is an
isomorphism of chain coalgebras.
Dualising provides $A \cong (UL)^{\dual}$ with the required $\Gamma$-algebra
structure.
\end{proof}

\bibliographystyle{amsplain}
\bibliography{homotopy}

\end{document}